\documentclass[12pt,twoside]{article}
\usepackage{graphicx}
\usepackage{amsfonts}
\usepackage{amsbsy}
\usepackage{amssymb}\usepackage{color}
\usepackage{amsmath}
\usepackage{cite}
\usepackage{pst-all}
\usepackage{pstricks}
\usepackage{graphics}
\def\bpsp{\begin{pspicture}}
\def\epsp{\end{pspicture}}

\newtheorem{theorem}{Theorem}[section]
\newtheorem{remark}[theorem]{Remark}
\newtheorem{example}[theorem]{Example}
\newtheorem{lemma}[theorem]{Lemma}
\newtheorem{corollary}[theorem]{Corollary}
\newtheorem{definition}[theorem]{Definition}
\newtheorem{proposition}[theorem]{Proposition}
\newtheorem{note}{Note}
\newtheorem{case}{Case}
\newtheorem{conjecture}{Conjecture}
\newtheorem{question}{Question}

\newcommand{\bea}{\begin{eqnarray}}
\newcommand{\eea}{\end{eqnarray}}
\newcommand{\beq}{\begin{eqnarray*}}
\newcommand{\eeq}{\end{eqnarray*}}

\def\m4{\mbox{\rm ~(mod $4$)}}

\def \bd{\begin{definition}}
\def \ed{\end{definition}}
\def \bqu{\begin{question}}
\def \equ{\end{question}}
\def \bcc{\begin{conjecture}}
\def \ecc{\end{conjecture}}
\def \bt{\begin{theorem}}
\def \et{\end{theorem}}
\def \bl{\begin{lemma}}
\def \el{\end{lemma}}
\def \bc{\begin{corollary}}
\def \ec{\end{corollary}}
\def \be{\begin{equation}}
\def \ee{\end{equation}}
\def \ben{\begin{enumerate}}
\def \een{\end{enumerate}}
\def \ba{\begin{array}}
\def \ea{\end{array}}
\def \bp{\begin{proposition}}
\def \ep{\end{proposition}}
\def \bx{\begin{example}}
\def \ex{\end{example}}
\def \br{\begin{remark}}
\def \er{\end{remark}}
\def \bdsc{\begin{description}}
\def \edsc{\end{description}}

\def \bn{\begin{case}}
\def \en{\end{case}}
\def \bnt{\begin{note}}
\def \ent{\end{note}}
\def\1{1\!\!1}

\def\mm2{\mbox{\rm ~(mod $2$)}}
\def\m4{\mbox{\rm ~(mod $4$)}}

\def\qed{\nolinebreak\hfill\rule{.2cm}{.2cm}\par\addvspace{.5cm}}

\def\m{\mu}

\def\1{\textbf{1}}
\def\0{\textbf{0}}

\begin{document}
\title{On $\alpha$-adjacency energy of graphs and Zagreb index}
\author{S. Pirzada$^a$, Bilal A. Rather$^b$, Hilal A. Ganie$^c$, Rezwan ul Shaban$^d$   \\
$^{a,b,c,d}${\em Department of Mathematics, University of Kashmir, Srinagar, India}\\
$^a$pirzadasd@kashmiruniversity.ac.in; ~~ $^b$bilalahmadrr@gmail.com\\
$^c$hilahmad1119kt@gmail.com, $^d$ rezwanbhat21@gmail.com}
\date{}

\pagestyle{myheadings} \markboth{}{On $\alpha$-adjacency energy of graphs}
\maketitle
\vskip 5mm
\noindent{\footnotesize \bf Abstract.} Let $A(G)$ be the adjacency matrix and $D(G)$ be the diagonal matrix of the vertex degrees of a simple connected graph $G$. Nikiforov defined the  matrix $A_{\alpha}(G)$ of the convex combinations of $D(G)$ and $A(G)$ as $A_{\alpha}(G)=\alpha D(G)+(1-\alpha)A(G)$,  for $0\leq \alpha\leq 1$.  If $ \rho_{1}\geq \rho_{2}\geq \dots \geq \rho_{n}$ are  the eigenvalues of $A_{\alpha}(G)$ (which we call $\alpha$-adjacency eigenvalues of $G$), the $ \alpha $-adjacency energy of $G$ is defined as $E^{A_{\alpha}}(G)=\sum_{i=1}^{n}\left|\rho_i-\frac{2\alpha m}{n}\right|$, where $n$ is the order and $m$ is the size of $G$. We obtain the upper and lower bounds for $E^{A_{\alpha}}(G) $ in terms of order $n$, size $m$ and Zagreb index $Zg(G)$ associated to the structure of $G$. Further, we characterize the extremal graphs attaining these bounds.
\vskip 3mm



\section{Introduction}

A simple graph is denoted by $G(V(G),E(G))$, where $V(G)=\{v_{1},v_{2},\ldots,v_{n}\}$ is its vertex set and ${E}(G)$ is its edge set. The \textit{order} and \textit{size} of $G$ are $|V(G)|=n$ and $|{E}(G)|=m$ respectively. The set of vertices adjacent to $v\in V(G)$, denoted by $N(v)$, refers to the \textit{neighborhood} of $v.$ The \textit{degree} of $v,$ denoted by $d_{G}(v)$ (we simply write $d_v$ if it is clear from the context) is the cardinality of $N(v)$. A graph is \textit{regular or degree regular} if all of its vertices are of the same degree. The adjacency matrix $A(G)=(a_{ij})$ of $G$ is a $(0, 1)$-square matrix of order $n$ whose $(i,j)$-entry is equal to 1, if $v_i$ is adjacent to $v_j$ and equal to 0, otherwise. If $D(G)={diag}(d_1, d_2, \dots, d_n)$ is the diagonal matrix of vertex degrees, the matrices $L(G)=D(G)-A(G)$ and $Q(G)=D(G)+A(G)$ are the Laplacian and the signless Laplacian matrices, respectively. Spectrum of $L(G)$ is the Laplacian spectrum and spectrum of $Q(G)$ is the signless Laplacian spectrum. The matrices $L(G)$ and $Q(G)$ are real symmetric and positive semi-definite. For $G$, we take $0=\mu_n\leq\mu_{n-1}\leq\dots\leq\mu_1$ and $0\leq q_n\leq q_{n-1}\leq\dots\leq q_1$ to be the Laplacian spectrum and signless Laplacin spectrum, respectively. For other standard notations, we refer to \cite{DMCM,hj,lsg,sp}.\\

Nikiforov \cite{vn} introduced the concept of \textit{merging $A$ and $Q$ spectral theories} by taking $A_{\alpha}(G)$ as the convex combinations of $D(G)$ and $A(G)$, and defined $A_{\alpha}(G)=\alpha D(G)+(1-\alpha)A(G)$,  for $0\leq \alpha\leq 1$. Since $A_{0}(G)=A(G), ~~~ 2A_{\frac{1}{2}}(G)=Q(G), ~~~ A_{1}(G)=D(G)$ and  $A_{\alpha}(G)-A_{\beta}(G)=(\alpha-\beta)L(G)$, any result regarding the spectral properties of $ A_\alpha $ matrix, has its counterpart for each of these particular graph matrices. Since the matrix $ A_{\alpha}(G)$ is real symmetric, all its eigenvalues are real and can be arranged  as $ \rho_{1}\geq \rho_{2}\geq \dots \geq \rho_{n}$. The largest eigenvalue  $ \rho_{1} $ (or simply $\rho(G)$) is called the \textit{spectral radius}. As $ A_{\alpha}(G) $ is nonnegative and irreducible, by the Perron-Frobenius theorem, $ \rho(G)$ is unique and there is a unique positive unit eigenvector $X$ corresponding to $ \rho(G),$ which is called the \textit{Perron vector} of $A_\alpha (G).$ Further results on spectral properties of the matrix  $A_\alpha (G)$ can be found in \cite{lxs,lhx,llx,ll,lds,ldss,vn,n01,xlls}.

Gutman \cite{IG1} defined the {\it energy} of a graph $G$ as $E(G)=\displaystyle\sum_{i=1}^{n}|\lambda_{i}|$, where $\lambda_{1}\geq\lambda_{2}\geq \dots\geq \lambda_{n}$ are the adjacency eigenvalues of $G$. Gutman et al. \cite{IGBZ} defined the {\it Laplacian energy} of a graph $G$ as $ LE(G)=\displaystyle\sum\limits_{i=1}^n\left|\mu_i-\frac{2 m)}{n}\right|$, $\mu_{1}\geq\mu_{2}\geq \dots\geq \mu_{n}$ are the Laplacian eigenvalues of $G$. For more details, see \cite{IG2}. Likewise, Abreu et al. \cite{abreu} defined the signless Laplacian energy of a graph $G$ as $QE(G)=\displaystyle\sum\limits_{i=1}^n\left|q_i-\frac{2 m}{n}\right|$, where $q_{1}\geq q_{2}\geq \dots\geq q_{n}$ are the signless Laplacian eigenvalues of $G$ and $\frac{2m}{n}$ is the average degree of $G$. For recent work, see \cite{hbp,ph}.\\
\indent Let $s_i=\rho_i-\frac{2\alpha m}{n}$ be the auxiliary eigenvalues corresponding to the eigenvalues of $ A_\alpha(G) $. The $\alpha$-adjacency energy $E^{A_{\alpha}}(G)$ \cite{hz} of a graph $G$ is defined as the mean deviation of the values of the eigenvalues of $A_{\alpha}(G)$, that is,
\begin{eqnarray}\label{1}
E^{A_{\alpha}}(G)=\sum_{i=1}^{n}\left|\rho_i-\frac{2\alpha m}{n}\right|=\displaystyle\sum_{i=1}^{n}|s_{i}|.
\end{eqnarray}
Obviously, $\sum_{i=1}^{n}s_{i}=0.$  From the definition, it is clear that $E^{A_{0}}(G)=E(G)$ and $2E^{A_{\frac{1}{2}}}(G)=QE(G)$.  Therefore, it follows that $\alpha$-adjacency energy of a graph $G$ merges the theories of (adjacency) energy and signless Laplacian energy. As such it will be interesting to study the quantity $E^{A_{\alpha}}(G)$. \\

The rest of the paper is organized as follows. In Section 2, we obtain  the upper bounds for $E^{A_{\alpha}}(G)$ and characterize the extremal graphs attaining these bounds. In Section 3,  we obtain  the lower bounds for $E^{A_{\alpha}}(G)$ and characterize the extremal graphs attaining these bounds.

\section{Upper bounds for $\alpha$-adjacency energy of a graph}

Let $\mathbb{M}_{m\times n}(\mathbb{R})$ be the set of all $m\times n$ matrices with real entries, that is, $\mathbb{M}_{m\times n}(\mathbb{R})=\{X: X=(x_{ij})_{m\times n}, x_{ij}\in \mathbb{R}\}$.
For $ M\in \mathbb{M}_{m\times n}(\mathbb{R}),$ the \emph{Frobenius norm} is defined as $$ \parallel M\parallel_F =\sqrt{\sum\limits_{i=1}^{n}\sum\limits_{j=1}^{n}|m_{ij}|^2}= \sqrt{trace(M^tM)} ,$$ where \emph{trace} of a square matrix is defined as sum of the diagonal entries. Further, if $ MM^t=M^tM $, then $ \parallel M\parallel_F^2=\sum\limits_{i=1}^{n}|\lambda_i(M)|^2 $, where $ \lambda_i $ is the $i^{th}$ eigenvalue of the matrix $ M$.\\
\indent The \emph{Zagreb index} $ Zg(G) $ of a graph $ G $ is defined as the sum of the squares of vertex degrees, that is, $ Zg(G)=\displaystyle \sum_{u\in V(G)} d_G^2(u) $.\\
\indent The following lemma can be found in \cite{hj}.
 \begin{lemma}\label{lem1}
  Let $X$ and $Y$ be Hermitian matrices of order $n$ and let $Z=X+Y$. Then
  \begin{align*}
  \lambda_k(Z)\leq \lambda_j(X)+\lambda_{k-j+1}(Y),~~ n\geq k\geq j\geq 1,\\
  \lambda_k(Z)\geq \lambda_j(X)+\lambda_{k-j+n}(Y), ~~n\geq j\geq k\geq 1,
  \end{align*}
 where $\lambda_i(M)$ is the $i^{th}$ largest eigenvalue of the matrix $M$. In either of these inequalities, equality holds if and only if there exists a unit vector which is an eigenvector corresponding to each of the three eigenvalues involved.\\
\end{lemma}

The following lemma gives some basic properties of the $\alpha$-adjacency matrix of $G$.

\begin{lemma}\label{lem2}
Let $G$ be a connected graph of order $n$ with $m$ edges and having vertex degrees $d_1\ge d_2\ge\dots\ge d_n$. Then\\
{\bf (1)}. $\sum\limits_{i=1}^{n}\rho_i=2\alpha m$ \quad {\bf (2)}. $\sum\limits_{i=1}^{n}\rho_i^2=\alpha^2 Zg(G)+(1-\alpha)^2\parallel A(G)\parallel_F^2$\\
{\bf (3)}. $\sum\limits_{i=1}^{n}s_i^2=\alpha^2Zg(G)+(1-\alpha)^2\parallel A(G)\parallel_F^2-\frac{4\alpha^2m^2}{n}$.\\
{\bf (4)}. $\rho(G)\geq \dfrac{2m}{n}$, equality holds if and only if $G$ is a degree regular graph.\\
{\bf (5)}. $ \rho(G)\geq \sqrt{\frac{Zg(G)}{n}} $, equality holds if and only if $G$ is a degree regular graph.
\end{lemma}
\textbf{Proof}. (1) Clearly, $\sum\limits_{i=1}^{n}\rho_i=\alpha \sum\limits_{i=1}^{n}d_i+(1-\alpha)\sum\limits_{i=1}^{n}\lambda_i=\alpha\displaystyle \sum_{u\in V(G)} d_G(u) =2\alpha m$.\\
(2). Here,
\begin{align*}
\sum\limits_{i=1}^{n}\rho_i^2& =\sum\limits_{i=1}^{n}\left(\alpha d_i+(1-\alpha)\lambda_i\right) ^2=\alpha^2\sum\limits_{i=1}^{n}d_i^2+(1-\alpha)^2\sum\limits_{i=1}^{n}(\lambda_i)^2 \\& =\alpha^2\displaystyle \sum_{u\in V(G)} d_G^2(u)+(1-\alpha)^22m =\alpha^2Z_g(G)+(1-\alpha)^2\parallel A(G)\parallel_F^2
\end{align*}
(3). We have,
\begin{align*}
\sum\limits_{i=1}^{n}s_i^2=\sum\limits_{i=1}^{n}\rho_i^2-\frac{4\alpha^2m^2}{n}
,\end{align*} and so by (2) the result follows.\\
(4). Let $\textbf{X}=\frac{1}{\sqrt{n}}(1,1,\dots,1) $ be a unit vector. Then, by Raleigh-Ritz's theorem for Hermitian matrices \cite{hj}, we have \[  \rho(G)\geq \dfrac{\textbf{X}^tA_\alpha(G)\textbf{X}}{\textbf{X}^t\textbf{X}}=\dfrac{\alpha\sum\limits_{i=1}^{n}d_i+
(1-\alpha)\sum\limits_{i=1}^{n}d_i}{n}=\dfrac{2m}{n}.\]
Assume that $G$ is $k$ degree regular. Then each row sum of $A_\alpha(G)$ equals to a constant $k$. Therefore, by the Perron-Frobenius theorem \cite{hj}, $k$ is a simple and largest eigenvalue of $A_\alpha(G)$. Thus $\rho(G)=k=\frac{nk}{n}=\frac{2m}{n}$ and equality holds. Conversely, assume that equality holds. Then $A_\alpha(G)\textbf{X} =\rho(G)\textbf{X}$. Therefore, $d_i=\rho(G)$ for all $i$ and thus $G$ is degree regular.\\
(5). This follows from \cite{vn}. Equality can be verified as in (4).\qed

\noindent From Case 3 of Lemma \ref{lem2}, we have $ \sum\limits_{i=1}^{n}s_i^2=(1-\alpha)^2\parallel A(G)\parallel_F^2+\sum\limits_{i=1}^{n}\left( \alpha d_i-\frac{2\alpha m}{n}\right)^2 $.\\
Let
\begin{equation}\label{bg}
2S(G)=(1-\alpha)^2\parallel A(G)\parallel_F^2+\sum\limits_{i=1}^{n}\left( \alpha d_i-\frac{2\alpha m}{n}\right) ^2.\end{equation}
We observe that $2S(G)=(1-\alpha)^2\parallel A(G)\parallel_F^2$ if and only if $ G $ is $\frac{2m}{n}$-degree regular graph, otherwise $ 2S(G)>(1-\alpha)^2\parallel A(G)\parallel_F^2 .$ Further $ 2S(G)=\parallel A(G)-\frac{2\alpha m}{n}\mathbb{I}_n\parallel_F^2=\sum\limits_{i=1}^{n}s_i^2 $, where $ \mathbb{I}_n $  is the identity matrix of order $ n .$ \\
\indent It is well known that a graph $G$ has two distinct eigenvalues if and only if $G\cong K_n$. Using this fact, it can be easily verified that the graph $G$ has two distinct $ \alpha $-adjacency eigenvalues if and only if $G$ is a complete graph with $ \alpha\ne 1 $. The $\alpha $-adjacency spectrum of the complete graph $ K_n $ is given in the next lemma \cite{vn}.

\begin{lemma}\label{lem5}
If $G=K_n$ is a complete graph, then the spectrum of $ A_\alpha(K_n)$  is $\{(n-1), (n\alpha-1)^{[n-1]}\} $, where $\rho^{[j]}$ means the eigenvalues $\rho$ is repeated $j$ times in the spectrum.
\end{lemma}

The following lemma \cite{vn} gives a lower bound for the $\alpha$-adjacency spectral radius.

\begin{lemma}\label{de}
If $ G $ is a graph with maximum degree $ \Delta(G)=\Delta $, then
\[ \rho(G)\geq\frac{1}{2}\left(\alpha(\Delta+1)+\sqrt{\alpha^2(\Delta+1)^2+4\Delta(1-2\alpha)} \right) .\]
For $ \alpha\in [0,1) $ and $ G $ being connected, equality holds if and only if $ G\cong K_{1,\Delta}. $
\end{lemma}

We first find the $\alpha$-adjacency energy of a degree regular graph.

\begin{theorem}\label{dreg}
If $ G $ is a degree regular graph of order $ n $ and $ \alpha\in\left [ 0,1\right )$, then
\[ E^{A_\alpha}(G)=(1-\alpha)E(G) .\]
\end{theorem}
\textbf{Proof:} Let $ \lambda_1, \lambda_2, \dots, \lambda_n $ be the adjacency eigenvalues of graph $ G .$ If $ G $ is a $k$ degree regular, then $ D(G)=kI_n $ and so \[ A_\alpha(G)=\alpha D(G)+(1-\alpha)A(G)=\alpha kI_n+(1-\alpha)A(G) .\] From this equality, it is clear that the $\alpha$-adjacency spectrum of $ G $ is $\{ \alpha k+(1-\alpha)\lambda_1, \dots,\alpha k+(1-\alpha)\lambda_n\}.$ Using this and the fact $\frac{2\alpha m}{n}=\alpha k$, we obtain $ E^{A_\alpha}(G)=(1-\alpha)E(G).$ \qed

From Theorem \ref{dreg}, for a degree regular graph $G$, it is clear that the value of $\alpha$-adjacency energy $E^{A_\alpha}(G)$ is a decreasing function of $\alpha$, for  $ \alpha\in\left [ 0,1\right )$.

 The following theorem gives McClelland type upper bound for $\alpha$-adjacency energy in terms of order $n$ and the quantity $S(G)$ associated to $G$.

\begin{theorem}\label{thm1}
If $ G $ is a connected graph of order $ n $, then $ E^{A_\alpha}(G)\leq \sqrt{2S(G)n}$.
\end{theorem}
\noindent{\bf Proof.} Using Cauchy-Schwarz's inequality, we have
\[ \left(E^{A_\alpha}(G)\right)^2=\left(\sum\limits_{i=1}^{n}|s_i|\right)^2\leq n\sum\limits_{i=1}^{n}s_i^2=2nS(G) \].\qed

Now, we obtain an upper bound for $\alpha$-adjacency energy in terms of order $n$, size $ m $ and the quantity $S(G)$ associated to $G$.

\begin{theorem}\label{thm2}
Let $ G $ be a connected graph of order $n\geq 3$ with $m$ edges and having Zagreb index $Zg(G)$. If  $\alpha\in [0,\frac{1}{2}]$ or $\alpha\in (\frac{1}{2},1)$ and $Zg(G)> \frac{8m^2}{n}-2m$ or $Zg(G)< \frac{4m^2}{n}$, then
\begin{equation}\label{b}
E^{A_\alpha}(G)\leq(1-\alpha)\left(  \dfrac{2m}{n}\right) +\sqrt{(n-1)\left[ 2S(G)-(1-\alpha)^2\left( \dfrac{2m}{n}\right)^2  \right] },
\end{equation}
where $ 2S(G) $ is same as in \eqref{bg}. Equality occurs if and only if either $ G= K_n $ or $ G $ is a connected  degree regular graph with three distinct eigenvalues given by $\frac{2m}{n} $, $\frac{2m\alpha}{n}+(1-\alpha) \sqrt{\frac{ 2m-\left( \frac{2m}{n}\right)^2}{n-1}} $ and $\frac{2m\alpha}{n}-(1-\alpha) \sqrt{\frac{ 2m-\left( \frac{2m}{n}\right)^2}{n-1}}$.
\end{theorem}
\textbf{Proof}. Let $\rho_1\ge \rho_2\ge\dots\ge \rho_n$ be $\alpha$-adjacency eigenvalues of $G$. For $1\le i\le n$, let $s_i=\rho_i(G)-\frac{2m\alpha}{n}$. Using Lemma \ref{lem2}, we have
$ \sum\limits_{i=2}^{n}s_i^2=2S(G)-s_1^2 .$ Applying Cauchy-Schwarz's inequality to the vectors $ (|s_2|, |s_3|,\dots,|s_n|) $ and $ (1,1,\dots,1) $, we obtain \[ \sum\limits_{i=2}^{n}|s_i|\leq\sqrt{ (n-1)\sum\limits_{i=2}^{n}s_i^2}=\sqrt{(n-1)\left[2S(G) -s_1^2 \right]}. \]
Therefore, we have
 \[E^{A_\alpha}(G)=s_1+\sum\limits_{i=2}^{n}|s_i|\leq s_1+ \sqrt{(n-1)\left[2S(G) -s_1^2 \right]}.\]
The last inequality suggests to consider the function $ F(x)= x+ \sqrt{(n-1)\left[2S(G) -x^2 \right] }.$ It is easy to see that this function is strictly decreasing in the interval $ \sqrt{2S(G)/n}<x\leq\sqrt{2S(G)} $. Since, $G$ is a connected graph, it follows that $m\ge n-1$ implying that  $ 2m\geq 2n-2> n $, for all $n\geq 3$.  We have $ \sqrt{2S(G)/n}\leq(1-\alpha)\frac{2m}{n}$ implying that
\begin{align}\label{tu}
\gamma \alpha^2-2\gamma^{'}\alpha+\gamma^{'}\geq 0,
\end{align} where $\gamma=\frac{8m^2}{n}-Zg(G)-2m$ and $\gamma^{'}=\frac{4m^2}{n}-2m$. For $\alpha=0$, inequality (\ref{tu}) follows, as $\gamma^{'}=\frac{4m^2}{n}-2m>0$. For $\alpha \in (0,1)$, consider the function $f(\alpha)=\gamma \alpha^2-2\gamma^{'}\alpha+\gamma^{'}$. It is easy to see that $f(\alpha)$ is decreasing for $\alpha\leq \frac{\gamma^{'}}{\gamma}$ and increasing for $\alpha\ge \frac{\gamma^{'}}{\gamma}$. If $Zg(G)> \frac{8m^2}{n}-2m$, then $\frac{\gamma^{'}}{\gamma}<0$, as $\gamma^{'}>0$ and so $\frac{\gamma^{'}}{\gamma}\notin (0,1)$. This gives $f(\alpha)> f(0)=\gamma^{'}>0$ and so inequality (\ref{tu}) follows in this case. So, assume that  $Zg(G)\le  \frac{8m^2}{n}-2m$. Then $\frac{\gamma^{'}}{\gamma}>0$. If $\frac{\gamma^{'}}{\gamma}\ge 1$, then $Zg(G)\geq \frac{4m^2}{n}$ and so it follows that $f(\alpha)\ge f(\frac{1}{2})=\frac{1}{4}\gamma >0$, for all $\frac{4m^2}{n}\le Zg(G)\le  \frac{8m^2}{n}-2m$. So, if $\frac{4m^2}{n}\le Zg(G)\le  \frac{8m^2}{n}-2m$, then inequality (\ref{tu}) holds for all $\alpha\in (0,\frac{1}{2}]$. Now, assume that $Zg(G)< \frac{4m^2}{n}$. It is clear that  $\frac{\gamma^{'}}{\gamma}\in (0,1)$ and so we have $f(\frac{\gamma^{'}}{\gamma})=\gamma^{'}\Big(1-\frac{\gamma^{'}}{\gamma}\Big)>0$, as $\frac{\gamma^{'}}{\gamma}<1$. So, if $Zg(G)< \frac{4m^2}{n}$, then inequality (\ref{tu}) holds for all $\alpha\in (0,1)$. Thus, it follows that the inequality  $\sqrt{2S(G)/n}\leq(1-\alpha)\frac{2m}{n}$ holds for all $\alpha\in [0,\frac{1}{2}]$ and holds for all $\alpha\in (\frac{1}{2},1)$, provided that $Zg(G)> \frac{8m^2}{n}-2m$ or $Zg(G)< \frac{4m^2}{n}$. Since $\rho_1\ge \frac{2m}{n}$, that is, $(1-\alpha)\frac{2m}{n}\leq s_1$, it follows that $ \sqrt{2S(G)/n}\leq(1-\alpha)\frac{2m}{n}\leq s_1 \leq \sqrt{2S(G)}$, for all $\alpha\in [0,\frac{1}{2}]$ and for all $\alpha\in (\frac{1}{2},1)$, provided that $Zg(G)> \frac{8m^2}{n}-2m$ or $Zg(G)< \frac{4m^2}{n}$. Now, $F(x)$ being decreasing in $ \sqrt{2S(G)/n}<x\leq\sqrt{2S(G)}$, it follows that  $ F(s_1)\leq F((1-\alpha)\frac{2m}{n}) $. Thus, from this, inequality \eqref{b} follows.\\
\indent  Suppose that equality occurs in \eqref{b}. Then all the inequalities above occur as equalities. By Lemma \ref{lem2}, equality occurs in $(1-\alpha)\frac{2m}{n}\leq s_1$, if and only if $G$ is a degree regular graph. Also,  equality occurs in Cauchy-Schwarz's inequality if $|s_2|=|s_3|=\cdots= |s_n|= \sqrt{\frac{ 2S(G)-(1-\alpha)^2\left( \frac{2m}{n}\right)^2 } {n-1}}$. Since, $\sqrt{2S(G)/n}\leq(1-\alpha)\frac{2m}{n}\le s_1$ holds for all $\alpha\in [0,\frac{1}{2}]$ and holds for all $\alpha\in (\frac{1}{2},1)$, provided that $Zg(G)> \frac{8m^2}{n}-2m$ or $Zg(G)< \frac{4m^2}{n}$, it follows that $ s_1> \sqrt{\frac{ 2S(G)-(1-\alpha)^2\left( \frac{2m}{n}\right)^2}{n-1}}$. Thus there are two cases  to consider. (i) Either $ G $ is a  connected degree regular graph with two distinct $\alpha$-adjacency eigenvalues (namely $\rho_1=\frac{2m}{n}$ and $\rho_2=\frac{2m\alpha}{n}-(1-\alpha)\sqrt{\frac{ 2m-\left( \frac{2m}{n}\right)^2}{n-1}}$ repeated $n-1$ times) or (ii) $ G $ is a connected degree regular graph with three distinct $\alpha$-adjacency  eigenvalues namely $\rho_1=\frac{2m}{n} $ and the other two given by $\frac{2m\alpha}{n}+(1-\alpha) \sqrt{\frac{ 2m-\left( \frac{2m}{n}\right)^2}{n-1}} $ and $\frac{2m\alpha}{n}-(1-\alpha) \sqrt{\frac{ 2m-\left( \frac{2m}{n}\right)^2}{n-1}} $. In Case (i), by Lemma \ref{lem5}, it follows that  $ G$ is a complete graph, that is  $G= K_n $, while in Case (ii),  it follows that $ G $ is a connected  degree regular graph with three distinct eigenvalues given by $\frac{2m}{n} $, $\frac{2m\alpha}{n}+(1-\alpha) \sqrt{\frac{ 2m-\left( \frac{2m}{n}\right)^2}{n-1}} $ and $\frac{2m\alpha}{n}-(1-\alpha) \sqrt{\frac{ 2m-\left( \frac{2m}{n}\right)^2}{n-1}}$.\\
\indent  Conversely, it can be easily verified that equality in \eqref{b} holds in each of above mentioned cases. \qed

Taking $\alpha=0$ and using the fact that $2S(G)=2m$,  we obtain the following result, which is the Koolen type \cite{k} upper bound for the energy $E(G)$.

\begin{corollary}\label{cor2}
Let $ G $ be a connected graph of order $ n\geq 3 $ with $m$ edges. Then
\begin{equation*}
E(G)\leq\frac{2m}{n} +\sqrt{(n-1)\left[2m-(\dfrac{2m}{n})^2  \right] }.
\end{equation*}
Equality occurs if and only if either $ G= K_n $ or $ G $ is a degree regular graph with  three distinct eigenvalues given by $ \frac{2m}{n}$ and other two with absolute value $ \sqrt{\frac{2m-\left(\frac{2m}{n}\right)^2}{n-1}}.$
\end{corollary}

Taking $\alpha=\frac{1}{2}$ and using the fact that $2S(G)=\frac{1}{4}\Big[2m+Zg(G)-\frac{4m^2}{n}\Big]$ together with $2E^{A_{\frac{1}{2}}}(G)=QE(G)$, we obtain the following result, which is the Koolen type upper bound for the signless Laplacian energy $QE(G)$.
\begin{corollary}\label{cor3}
Let $ G $ be a connected graph of order $ n \ge 3$ with $m$ edges having Zagreb index $Zg(G)$. Then
\begin{equation*}
QE(G)\leq\frac{2m}{n} +\sqrt{(n-1)\left[2m+Zg(G)-\frac{4m^2}{n}\left(1+\frac{1}{n}\right) \right]}.
\end{equation*}
Equality occurs if and only if either $ G= K_n $ or $ G $ is a degree regular graph with  three distinct eigenvalues given by $ \frac{4m}{n}$, $\frac{2m}{n}+ \sqrt{\frac{ 2m-\left( \frac{2m}{n}\right)^2}{n-1}} $ and $\frac{2m}{n}+ \sqrt{\frac{ 2m-\left( \frac{2m}{n}\right)^2}{n-1}} $.
\end{corollary}

The following lemma gives a relation between $\alpha$-adjacency eigenvalues of $G$ and $\alpha$-adjacency eigenvalues of  spanning subgraphs of $G$.

\begin{lemma}\label{lem3}
Let $ G $ be a connected graph of order $ n\ge 3$ and let $ \alpha \in [\frac{1}{2},1)$. If $ G^{'} $ is the graph obtained from $ G $ by deleting an edge, then for any $ 1\leq i\leq n$, we have $ \rho_i(G)\geq \rho_i (G^{'}).$
\end{lemma}
\textbf{Proof}. Let $G$ be a connected graph of order $n\geq 3$ and let $e=uv$ be an edge in $G$. Let $G^{'}=G-e$ be the graph obtained from $G$ by deleting $e$. It is easy to see that
\begin{align}\label{t}
A_\alpha(G)=A_\alpha(G^{'})+N,
\end{align}
 where $N$ is the matrix of order $n$ indexed by the vertices of $G$ having $(u,v)^{th}$ and $(v,u)^{th}$ entries both equal to $1-\alpha$, and the $(u,u)^{th}$ and $(v,v)^{th}$ entries both equal to $\alpha$, and all other entries equal to zero. It can be seen that the eigenvalues of the matrix $N$ are $1^{[1]},2\alpha-1^{[1]},0^{[n-2]},$ where $\lambda^{[j]}$ means the eigenvalue $\lambda$ is repeated $j$ times in the spectrum. Taking $Z=A_\alpha(G), X=A_\alpha(G^{'}), Y=N$ and $k=j=i$ in the second inequality of Lemma \ref{lem1}, we get $ \rho_i(G)\geq \rho_i (G^{'}),$ provided that $ \alpha \in [\frac{1}{2},1)$.\qed

In $G$, let $\eta=\eta(G)$ be the number of $\alpha$-adjacency eigenvalues greater or equal to $\frac{2m\alpha}{n}$. Since, by Lemma \ref{lem2}, we have $\rho_1\ge \frac{2m}{n}$, it follows that $1\le \eta\le n$. Parameters similar to $\eta$ have been considered for the graph matrices and therefore it will be interesting to connect the parameter $\eta$ with $\alpha$-adjacency energy of $G$. Now, we obtain an upper bound for $\alpha$-adjacency energy in terms of order $n$, size $m$ and the parameter $\eta$  associated to $G$.
\begin{theorem}\label{thm31}
Let $ G $ be a connected graph of order $ n\ge 3 $ and  let  $ \frac{1}{2}\leq \alpha<  1 .$ Then
\begin{align*}
E^{A_\alpha}(G)\leq 2(n-1) +2(\eta-1)(\alpha n-1)-\frac{4\alpha \eta m}{n},
\end{align*}
with equality if and only if $ G\cong K_n .$
\end{theorem}
\textbf{Proof.} Let $ G $ be a connected graph of order $ n $ having $\alpha$-adjacency eigenvalues $ \rho_1\geq\rho_2\geq\dots\geq\rho_n .$ Let $ \eta $ be the positive integer such that $ \rho_{\eta}\geq\frac{2\alpha m}{n}$ and $ \rho_{\eta+1}<\frac{2\alpha m}{n} $. Using (1) of Lemma \ref{lem2} and the definition of $\alpha$-adjacency energy, we have
\begin{align*}
E^{A_\alpha}(G)&=\sum\limits_{i=1}^{n}\left|\rho_i-\frac{2\alpha m}{n}\right|=\sum\limits_{i=1}^{\eta}\left( \rho_i-\frac{2\alpha m}{n}\right) +\sum\limits_{i=\eta+1}^{n}\left( \frac{2\alpha m}{n}-\rho_i\right)\\
&=2\left(\sum\limits_{i=1}^{\eta}\rho_i-\frac{2\eta\alpha m}{n} \right).
\end{align*}
 Clearly $G$ is a spanning subgraph of $ K_n$. So from Lemma \ref{lem3}, it follows that $ \rho_i(G)\leq \rho_i(K_n) $ for each $ 1\leq i\leq n $. Therefore, we have
\begin{equation}\label{mm}
\sum\limits_{i=1}^{\eta}\rho_i(G)\leq \sum\limits_{i=1}^{\eta}\rho_i(K_n)=n-1+(\eta-1)(\alpha n-1).
\end{equation}
Using \eqref{mm}, we obtain
\begin{align*}
E^{A_\alpha}(G)\leq 2(n-1) +2(\eta-1)(\alpha n-1)-\frac{4\alpha \eta m}{n}.
\end{align*}
 Assume that equality occurs so that equality occurs in \eqref{mm}. Since,  equality occurs  in \eqref{mm} if and only if $ G\cong K_n $, it follows that equality holds if and only if $ G\cong K_n.$ \qed

The following lemma will be required in the sequel.

\begin{lemma}\label{thul1}
Let $G$ be a connected graph of order $n$, size $m$ and having vertex degrees $d_1\ge d_2\ge \dots\ge d_n$. Then
\[  \rho(G)\geq \sqrt{\frac{Zg(G)}{n}} \geq \frac{2m}{n}.\]
\end{lemma}
\textbf{Proof.} The first inequality follows by Case 5 of Lemma \ref{lem2}. Therefore, we need to prove the  second inequality. Applying Cauchy-Schwartz inequality to $\left( \sum\limits_{i=1}^{n}d_i \right)^2$, we have $ \left( \sum\limits_{i=1}^{n}d_i \right)^2\leq n\sum\limits_{i=1}^{n}d_i^2 $, which implies $ \sum\limits_{i=1}^{n}d_i^2\geq \frac{(2m)^2}{n} $ and hence $ \sqrt{\sum\limits_{i=1}^{n}d_i^2}\geq\frac{2m}{\sqrt{n}} .$ Thus,
\[ \sqrt{\dfrac{Zg(G)}{n}}= \sqrt{\dfrac{\sum\limits_{i=1}^{n}d_i^2}{n}} \geq \dfrac{2m}{n} .\] This completes the proof.\qed

The following theorems give upper bounds for the $\alpha$-adjacency energy in terms of order $n$, size $m$, Zagreb index $Zg(G) $ and the parameter $\alpha$.

\begin{theorem}\label{thu2}
Let $ G $ be a connected graph of order $ n\geq 3 $ having Zagreb index $Zg(G)$ and let $ \alpha\leq 1-\frac{n}{2m}.$ Then
\begin{equation}\label{eq4.7}\begin{split}
E^{A_\alpha}(G) &\leq  \alpha^2Zg(G)+(1-\alpha )^2\parallel A(G)\parallel_F^2-\frac{2\alpha m}{n^2}(2\alpha nm+2\alpha m+n)+ \ln\left( \frac{\theta}{\Gamma}\right)\\
&\quad+\frac{4\alpha m}{n}\left(\sqrt{\frac{Zg(G)}{n}} \right) -\left( \sqrt{\frac{Zg(G)}{n}}\right) \left( \sqrt{\frac{Zg(G)}{n}}-1\right),
\end{split}
\end{equation}
where $ \Gamma=\left| det\left(A_\alpha(G)-\frac{2\alpha m}{n}\mathbb{I}_n \right) \right| $ and $ \theta=\sqrt{\frac{Zg(G)}{n}}-\frac{2\alpha m}{n} .$  Equality holds if and only if $ G \cong K_n$ and $\alpha=0$ or $ G $ is a $ k$-degree regular graph with three distinct $\alpha$-adjacency eigenvalues given by $k, \alpha k+1$ and $\alpha k-1$.
\end{theorem}
\noindent\textbf{Proof.} Let $ G $ be a connected graph of order $ n$ and let $ \rho_1\geq\rho_2\geq\dots\geq\rho_n $ be the $ \alpha $-adjacency eigenvalues of $ G $.
Consider the function
\[f(x)=\left(x-\frac{2\alpha m}{n}\right)^2-\left(x-\frac{2\alpha m}{n}\right)-\ln\left(x-\frac{2\alpha m}{n}\right),\ \ \ \  \left(x-\frac{2\alpha m}{n}\right)>0.\]
It is easy to see that this function is non-decreasing for $ x-\frac{2\alpha m}{n}\geq 1 $ and non-increasing for $ 0\leq \left(x-\frac{2\alpha m}{n}\right)\leq 1 $. So, we have  $ f(x)\geq f\left(\frac{2\alpha m}{n}+1\right)=0 $ implying that
$$ \left(x-\frac{2\alpha m}{n}\right)\leq \left(x-\frac{2\alpha m}{n}\right)^2-\ln\left(x-\frac{2\alpha m}{n}\right) $$
for $ \left(x-\frac{2\alpha m}{n}\right)>0, $ with equality if and only if $ \left(x-\frac{2\alpha m}{n}\right)=1. $ Using  these observations in the  definition of $ \alpha $-adjacency energy, we have
\begin{align}\label{up1}
E^{A_\alpha}(G)&=\rho_1-\frac{2\alpha m}{n}+\sum\limits_{i=2}^{n}\left|\rho_i-\frac{2\alpha m}{n}\right|\nonumber\\
&\leq\rho_1-\frac{2\alpha m}{n}+\sum\limits_{i=2}^{n}\left( \left( \rho_i-\frac{2\alpha m}{n}\right)^2-\ln\left|\rho_i-\frac{2\alpha m}{n}\right|\right)\nonumber\\&
= \rho_1-\frac{2\alpha m}{n}+\alpha^2 Zg(G)+(1-\alpha )^2\parallel A(G)\parallel_F^2- \left(\frac{4\alpha^2 m^2}{n^2}\right)(n-1)-\rho_1^2\nonumber\\
&\quad-\ln\prod\limits_{i=1}^{n}\left| \rho_i-\frac{2\alpha m}{n}\right|+\ln\left(\rho_1-\frac{2\alpha m}{n}\right) -\frac{4\alpha m}{n}(2\alpha m-\rho_1)\nonumber\\
&= \alpha^2 Zg(G)+(1-\alpha )^2\parallel A(G)\parallel_F^2-\frac{2\alpha m}{n^2}(2\alpha nm+2\alpha m+n)+\frac{4\alpha m}{n}\rho_1\nonumber\\
&\quad- \ln\Gamma +\ln\left(\rho_1-\frac{2\alpha m}{n}\right)-\rho_1(\rho_1-1).
\end{align}
Consider the function \[\begin{split}
g(x)&=\alpha^2Zg(G)+(1-\alpha )^2\parallel A(G)\parallel_F^2-\frac{2\alpha m}{n^2}(2\alpha nm+2\alpha m+n)+\frac{4\alpha m}{n}x\\
&\quad- \ln\Gamma +\ln\left(x-\frac{2\alpha m}{n}\right)-x(x-1).
\end{split}\]
Evidently, the function $g(x)$ is increasing for $ 0\leq x-\frac{2\alpha m}{n}\leq 1 $ and decreasing for $ x-\frac{2\alpha m}{n}\geq 1 $. Since  $ x-\frac{2\alpha m}{n}\geq(1-\alpha)\frac{2m}{n}\geq 1 $ provided that $\alpha\le 1-\frac{n}{2m}$, then t for $\alpha\le 1-\frac{n}{2m}$, it follows that $ x-\frac{2\alpha m}{n}\ge 1$. Further, $(1-\alpha)\frac{2m}{n}\geq 1$ implies that $\frac{2m}{n}\geq 1+\frac{2m\alpha}{n}$ and by Lemma \ref{thul1}, we have $x\ge \sqrt{\frac{Zg(G)}{n}} \geq \dfrac{2m}{n}$. Therefore, it follows that
\begin{align}\label{up2}
g(x)&\leq g\left(\sqrt{\frac{Zg(G)}{n}}\right) = \alpha^2 Zg(G)+(1-\alpha )^2\parallel A(G)\parallel_F^2-\frac{2\alpha m}{n^2}(2\alpha nm+2\alpha m+n)\nonumber\\
&\quad+\frac{4\alpha m}{n}\sqrt{\frac{Zg(G)}{n}} + \ln\left( \frac{\theta}{\Gamma}\right)- \frac{Zg(G)}{n}-\sqrt{\frac{Zg(G)}{n}}.
\end{align}
Combining inequalities \eqref{up2} and \eqref{up1}, we arrive at \eqref{eq4.7}.\\
\indent Assume that inequality holds in \eqref{thu2}. Then all the inequalities above occur as equalities. By Lemma \ref{lem2}, equality  occurs in $ \rho_1\ge \sqrt{\frac{Zg(G)}{n}}$ if and only if $G$  is a degree regular graph. For equality in \eqref{up1}, we have $ \left| \rho_2-\frac{2\alpha m}{n} \right|=\dots=\left| \rho_n-\frac{2\alpha m}{n} \right| =1 $. For $i=2,3,\dots,n$, the quantity  $ \left|\rho_{i}-\frac{2\alpha m}{n}\right| $ can have at most two distinct values and therefore we have the following cases.\\
\textbf{Case 1.} For all $i=2,3,\dots,n$, if $ \rho_{i}-\frac{2\alpha m}{n}=1$, then $\rho_i=1+\frac{2\alpha m}{n}$, implying that $G$ has two distinct $\alpha$-adjacency  eigenvalues, namely $\rho_1=\frac{2m}{n}$ and $\rho_i=1+\frac{2\alpha m}{n}$. So, by Lemma \ref{lem5},  equality occurs for the complete graph $K_n$, provided that $\alpha$-adjacency eigenvalues of $K_n$ are $n-1$ with multiplicity $1$ and $\alpha n-1$ with multiplicity $n-1$. It is clear  that equality can  not hold in this case.\\
\textbf{Case 2.} For all $i=2,3,\dots,n$, if $ \rho_{i}-\frac{2\alpha m}{n}=-1$, then $\rho_i=\frac{2\alpha m}{n}-1$, implying that $G$ has two distinct $\alpha$-adjacency  eigenvalues, namely $\rho_1=\frac{2m}{n}$ and $\rho_i=\frac{2\alpha m}{n}-1$. So, using Lemma \ref{lem5}, it follows that equality occurs for the complete graph $K_n$, provided that $\alpha=0$.\\
{\bf Case 3.} For the remaining case, for some $t$, let $ \rho_{i}-\frac{2\alpha m}{n}=1,$ for $i=2,3,\dots,t$, and $ \rho_{i}-\frac{2\alpha m}{n}=-1$, for $i=t+1,\dots,n$. This implies that $G$ is degree regular graph with three distinct $\alpha$-adjacency eigenvalues, namely $\rho_1=\frac{2 m}{n}$ with multiplicity $1$, $\rho_{i}=1+\alpha \rho_1$ with multiplicity $t-1$ and $\rho_{i}=\alpha \rho_1-1$ with multiplicity $n-t$.\\
\indent Conversely, if $G\cong K_n$, then $\rho_1=n-1$, $\rho_i=\alpha n-1$, for $i=2,3,\dots,n$ and $\frac{2\alpha m}{n}=\alpha(n-1)$. It can be seen that equality occurs in \eqref{thu2}. On the other hand, if $G$ is a degree regular graph with three distinct $\alpha$-adjacency eigenvalues, namely $\rho_1,~ \alpha \rho_1+1$ and $\alpha \rho_1-1$, then from the above discussion, it is clear that the equality holds in \eqref{thu2}. \qed

\begin{theorem}\label{thu1}
Let $ G $ be a connected graph of order $ n\geq 3 $ having Zagreb index $Zg(G)$  and let $ \alpha\leq 1-\frac{n}{2m}.$ Then
\begin{equation*}\begin{split}
E^{A_\alpha}(G) &\leq \alpha^2Zg(G)+(1-\alpha)^2\parallel A(G)\parallel_F^2 +\ln\left( \frac{2m(1-\alpha)}{n\Gamma} \right)\\
&\quad -\frac{2\alpha m}{n^2}\left( 2n\alpha m+2\alpha m-4m+n \right)-\frac{2m}{n^2}\left( 2m-n \right),
\end{split}
\end{equation*}
where $ \Gamma=\left| det\left(A_\alpha(G)-\frac{2\alpha m}{n}\mathbb{I}_n \right) \right|.$ Equality holds if and only if $ G \cong K_n$ and $\alpha=0$ or $ G $ is a $ k$-degree regular graph with three distinct $\alpha$-adjacency eigenvalues given by $k,~ \alpha k+1$ and $\alpha k-1$.
\end{theorem}
\noindent{\bf Proof.} The proof is similar to the proof of Theorem \ref{thu2}.

\section{Lower bounds for the $\alpha$-adjacency energy of graphs}

The following theorem gives a lower bound for the $\alpha$-adjacency energy.

\begin{theorem}
If $ G $ is a connected graph of order $ n\geq 3 $, size $m$ and Zagreb index $Zg(G)$, then
\begin{equation}\label{lbe}
E^{A_\alpha}(G)\geq \sqrt{2\left( \alpha^2Zg(G)+(1-\alpha)^2\parallel A(G)\parallel_F^2-\frac{2(\alpha m)^2}{n}\right) }.
\end{equation}
for  $ \alpha\in [0,1) $
\end{theorem}
\noindent\textbf{Proof.} Let $s_1\ge s_2\ge\dots\ge s_n$ be the auxiliary eigenvalues as defined earlier. We have
\[ \left( E^{A_\alpha}(G)\right)^2 =\left( \sum\limits_{i=1}^{n}|s_i|\right)^2=\sum\limits_{i=n}^{n}s_i^2(G) +2\sum\limits_{1\leq i<j\leq n}|s_is_j| . \] By Case 3 of  Lemma \ref{lem2}, we have
\begin{equation}\label{lbe1}
\sum\limits_{i=1}^{n}s_i^2(G)=\alpha^2Zg(G)+(1-\alpha)^2\parallel A(G)\parallel_F^2-\dfrac{4\alpha^2m^2}{n}.
\end{equation} Also,
\begin{align*}
&2\sum\limits_{1\leq i<j\leq n}|s_is_j|
\geq 2\left|\sum\limits_{1\leq i<j\leq n}\left( \rho_i(G)-\frac{2\alpha m}{n}\right)\left( \rho_j(G)-\frac{2\alpha m}{n}\right)\right|\\
& =2\left|\sum\limits_{1\leq i<j\leq n}\rho_i(G)\rho_j(G)-\frac{4\alpha^2m^2}{n}(n-1) +\frac{4\alpha^2m^2}{n^2}n(n-1)\right|\\
& = 2\left|\sum\limits_{i=1}^{n}\rho_i(G)\rho_j(G)\right|=\sum\limits_{i=1}^{n}\rho_i^2(G) =\alpha^2\sum\limits_{i=1}^{n}d_i^2+(1-\alpha)^2\parallel A(G) \parallel_F^2.
\end{align*}
Using this inequality and \eqref{lbe1}, clearly \eqref{lbe} follows.
\qed

Now, we obtain a lower bound for the $\alpha$-adjacency energy in terms of order $n$, size $m$ and the parameter $\alpha$.

\begin{theorem}\label{thm41}
Let $ G $ be a connected graph of order $ n \ge 3$ and size $m$ and let $ \alpha\in [0,1) $. Then
\begin{equation}\label{lbe2}
E^{A_\alpha}(G)\geq 4(1-\alpha)\frac{m}{n}.
\end{equation}
Equality occurs if and only if $ G $ is degree regular with one positive and $n-1$ negative adjacency eigenvalues.
\end{theorem}
\textbf{Proof.} Let $ G $ be a connected graph of order $n$ and having $\alpha$-adjacency eigenvalues $ \rho_1\geq\rho_2\geq\cdots\geq\rho_n $. Let $ \eta $ be a positive integer such that $ \rho_{\eta}\geq\frac{2\alpha m}{n}$ and $ \rho_{\eta+1}<\frac{2\alpha m}{n} $. Using Case 1 of Lemma \ref{lem2} and the definition of $\alpha$-adjacency energy, we have
\begin{align*}
E^{A_\alpha}(G)&=\sum\limits_{i=1}^{n}\left|\rho_i-\frac{2\alpha W(G)}{n}\right|=\sum\limits_{i=1}^{\eta}\left( \rho_i-\frac{2\alpha m}{n}\right) +\sum\limits_{i=\eta+1}^{n}\left( \frac{2\alpha m}{n}-\rho_i\right)\\
&=2\left(\sum\limits_{i=1}^{\eta}\rho_i-\frac{2\eta\alpha m}{n} \right).
\end{align*}
First we show that
\begin{equation}\label{def}
E^{A_\alpha}(G)=2\left(\sum\limits_{i=1}^{\eta}\rho_i-\frac{2\eta\alpha m}{n} \right)=2 \max_{1\leq j\leq n}\left\lbrace \sum\limits_{i=1}^{j}\rho_i-\frac{2\alpha jm }{n}\right\rbrace.
\end{equation}
Since $1\le \eta\leq n$, it follows that either $\eta<j$ or $\eta\ge j$. If  $ j>\eta$, then we have
\begin{align*}
\sum\limits_{i=1}^{j}\rho_i -\frac{2\alpha jm}{n} &
=\sum\limits_{i=1}^{\eta}\rho_i+\sum\limits_{i=\eta+1}^{j}\rho_i-\frac{2\alpha jm}{n} \\
& < \sum\limits_{i=1}^{\eta}\rho_i-\frac{2\alpha \eta m}{n}
\end{align*}
as $\rho_i<\frac{2\alpha m}{n}$, for $i\geq \eta+1$. Similarly, for $ j\leq \eta$, it can be seen that
\begin{align*}
\sum\limits_{i=1}^{j}\rho_i -\frac{2\alpha jm}{n}\leq \sum\limits_{i=1}^{\eta}\rho_i-\frac{2\alpha \eta m}{n}.
\end{align*}
This proves \eqref{def}. Therefore, we have
\begin{align*}
E^{A_\alpha}(G)&= 2\max_{1\leq j\leq n}\left\lbrace \sum\limits_{i=1}^{j}\rho_i-\frac{4\alpha jm }{n}\right\rbrace\geq 2\rho_1-\frac{4\alpha m}{n}\\
&\geq\frac{4 m}{n}-\frac{4\alpha m}{n}=(1-\alpha)\frac{4 m}{n}.
\end{align*}
Suppose equality holds in \eqref{lbe2}. Then all the inequalities above occur as equalities. By Lemma \ref{lem2}, equality occurs in $\rho_1\ge \frac{2m}{n}$  if and only if $G$ is a degree regular graph. Also, equality occurs in $2\max_{1\leq j\leq n}\left\lbrace \sum\limits_{i=1}^{j}\rho_i-\frac{4\alpha jm }{n}\right\rbrace\geq 2\rho_1-\frac{4\alpha m}{n}$ if and only if $\eta=1$. Thus, it follows that equality occurs in \eqref{lbe2} if and only if $G$ is a degree regular graph with $\eta=1$. Let $G$ be a $k$-degree regular graph having adjacency eigenvalues $\lambda_1\geq \lambda_2\geq \dots\geq \lambda_n$. Then, by Theorem \ref{dreg}, we have $\rho_1=k$ and $\rho_i=\alpha k+(1-\alpha)\lambda_i$, for $i=2,3,\dots,n$. Since $\eta=1$, for $2\le i\le n$, we have  $ \rho_{i}<\frac{2\alpha m}{n}=\alpha k$, which gives $\alpha k+(1-\alpha)\lambda_i<\alpha k$, which further gives $\lambda_i<0$ as $1-\alpha>0$. Thus, it follows that equality  occurs in \eqref{lbe2} if and only if $G$ is a degree regular graph with one positive and $n-1$ negative adjacency eigenvalues.\qed

Proceeding similarly as in Theorem  \ref{thm41} and using Case 5 of Lemma \ref{lem2}, we obtain the following lower bound for  $\alpha$-adjacency energy of a connected graph.

\begin{theorem}
Let $ G $ be a connected graph of order $ n\ge 3 $, size $m$  and Zagreb index $Zg(G)$. For $ \alpha\in [0,1) $, we have
\begin{equation}
E^{A_\alpha}(G)\geq 2\sqrt{\frac{Zg(G)}{n}}-\frac{4\alpha m}{n}.
\end{equation}
Equality occurs if and only if $ G $ is degree regular with one positive and $n-1$ negative adjacency eigenvalues.
\end{theorem}

The following theorem gives a lower bound for $\alpha$-adjacency energy of a connected graph in terms of order $n$, size $m$, the maximum degree $\Delta$ and the parameter $\alpha$.
\begin{theorem}
Let $ G $ be a connected graph of order $ n\ge 3 $, size $m$ and maximum degree $\Delta.$ For $ \alpha\in [0,1) $, we have
\begin{equation}\label{delb}
E^{A_\alpha}(G)\geq \left(\alpha(\bigtriangleup+1)+\sqrt{\alpha^2(\bigtriangleup+1)^2+4\bigtriangleup(1-2\alpha)} \right) -\frac{4\alpha m}{n},\end{equation}
with equality if and only if $ G\cong K_{1,\Delta} .$
\end{theorem}
\textbf{Proof.} By Equation \eqref{def} and Lemma \ref{de}, we have
\begin{align*}
E^{A_\alpha}(G)&= \max_{1\leq j\leq n}\left\lbrace 2\sum\limits_{i=1}^{j}\rho_i-\frac{4\alpha im }{n}\right\rbrace\geq 2\rho_1(G)-\frac{4\alpha m}{n}\\
&\geq\alpha(\Delta+1)+\sqrt{\alpha^2(\Delta+1)^2+4\Delta(1-2\alpha)} -\frac{4\alpha m}{n}.
\end{align*}
Suppose equality holds in \eqref{delb}. Then all the inequalities above occur as equalities. Since  equality occurs in Lemma \ref{de} if and only if $G\cong K_{1,\Delta}$ and equality occurs in
$$ \max_{1\leq j\leq n}\left\lbrace 2\sum\limits_{i=1}^{j}\rho_i-\frac{4\alpha im }{n}\right\rbrace\geq 2\rho_1(G)-\frac{4\alpha m}{n}$$
if and only if $\eta=1$, it follows that equality occurs in \eqref{delb} if and only if $G\cong K_{1,\Delta}$, $\Delta=n-1$  and $\eta=1$. For the graph $ K_{1,\Delta}$, the $ \alpha $-adjacency eigenvalues are $\left\lbrace  \alpha^{[n-2]}, \frac{\alpha (\Delta+1)\pm\sqrt{D}}{2} \right\rbrace $, where $ D=\alpha^2(\Delta+1)^2+4\Delta(1-2\alpha) $ and average of the $\alpha$-adjacency  equals to $ 2\alpha -\frac{2\alpha}{n}  $. Clearly, now $ \eta=1 $ for $ K_{1,\Delta} .$ Thus equality occurs in \eqref{lbe2} if and if $ G\cong K_{1,\Delta} .$ This completes the proof. \qed

Now, we obtain a lower bound for $\alpha$-adjacency energy of a connected graph in terms of order $n$, size $m$ and Zagreb index $Zg(G)$.

\begin{theorem}
Let  $G$ be a connected graph of order $n\ge 3$ having size $m$ and Zagreb index $Zg(G)$. For $ \alpha\in [0,1) $, we have
\begin{equation}\label{lp1}
E^{A_{\alpha}}(G)\geq \sqrt{\frac{Zg(G)}{n}}+(n-1)+\ln \left( \frac{\Gamma}{\theta} \right),
\end{equation}
where $ \Gamma=\left|det\left( A_\alpha(G)-\frac{2\alpha m}{n}\mathbb{I}_n\right) \right| $ and $ \theta=\sqrt{\frac{Zg(G)}{n}}-\frac{2m\alpha}{n} $. Equality holds as in Theorem \ref{thu2}.
\end{theorem}
\noindent\textbf{Proof.} Consider the function $ f(x)=x-1-\ln x $, where $ x>0 $. It is easy to verify that the function $ f(x) $ is increasing for $ x\geq 1 $ and decreasing for $ 0\leq x\leq 1 $. Therefore, we have $ f(x)\geq f(1)=0 $ implying that $ x\geq 1+\ln x $, for $ x>0 $, with equality  if and only if $ x=1.$ Using this observation with $x=\Big|\rho_i-\frac{2\alpha m}{n}\Big|$, for $2\le i\le n$ and the definition of $\alpha$-adjacency energy, we have
\begin{align*}
E^{A_\alpha}(G)&=\rho_1-\frac{2\alpha m}{n}+\sum\limits_{i=2}^{n}\left|\rho_i-\frac{2\alpha m}{n}\right|\\
&\geq\rho_1-\frac{2\alpha m}{n}+(n-1)+\sum\limits_{i=2}^{n}\ln\left|\rho_i-\frac{2\alpha m}{n}\right| \\
&= \rho_1-\frac{2\alpha m}{n}+(n-1)+\ln\prod\limits_{i=2}^{n}\left|\rho_i-\frac{2\alpha m}{n}\right|\\
&= \rho_1-\frac{2\alpha m}{n}+(n-1)+\ln\left|det\left(A_\alpha(G)-\frac{2\alpha m}{n}\right) \right|-\ln\left( \rho_1-\frac{2\alpha m}{n}\right).
\end{align*}
Now, consider the function $ g(x)=x-\frac{2\alpha m}{n}+(n-1)+\ln\left|det\left(A_\alpha(G)-\frac{2\alpha m}{n}\right) \right|-\ln\left( x-\frac{2\alpha m}{n}\right) $. Clearly, $ g(x) $ is increasing for $x-\frac{2\alpha m}{n}\ge 1$. Since, $x-\frac{2\alpha m}{n}\ge (1-\alpha)\frac{2m}{n}\ge 1$ implying that $\alpha\leq 1-\frac{n}{2m}$, therefore, for $\alpha\leq 1-\frac{n}{2m}$, it follows that $x-\frac{2\alpha m}{n}\ge 1$. Further, $(1-\alpha)\frac{2m}{n}\ge 1$ implies that $\frac{2m}{n}\ge 1+\frac{2m\alpha}{n}$. From Lemma \ref{lem2} and the fact that $g(x)$ is increasing for $1+\frac{2m\alpha}{n}\le\frac{2m}{n}\le \sqrt{\frac{Zg(G)}{n}}\le x$, it follows that $ g(x)\geq g(\sqrt{\frac{Zg(G)}{n}})$. From this, Inequality \eqref{lp1} follows. Equality case can be discussed similar to Theorem \ref{thu2}.\qed

A lower bound similar to the lower bound given in the above theorem can be obtained for $ \rho\geq \frac{2m}{n}$.\\

\noindent{\bf Acknowledgements.} The research of S. Pirzada is supported by SERB-DST, New Delhi under the research project number  MTR/2017/000084.

\end{document}